\documentclass[11pt,a4paper]{article}
\usepackage[pdftex]{graphicx}
\usepackage{amsmath,amssymb,amsfonts,amsthm}
\numberwithin{equation}{section}
\usepackage{indentfirst}
\usepackage{enumitem} % for customizing lists
\usepackage[margin=1.3in]{geometry}
\usepackage{fancyhdr}

\usepackage[colorlinks=true,linkcolor=magenta,citecolor=blue,urlcolor=cyan]{hyperref}
\usepackage{titlesec}
\usepackage{etoolbox}
\usepackage{graphicx}  % for scalebox if needed later
\usepackage{changepage}
\usepackage{float}
\theoremstyle{plain}
\newtheorem{theorem}{Theorem}[section]
\newtheorem{lemma}[theorem]{Lemma}

\usepackage{comment}

\makeatletter
\renewcommand{\maketitle}{
	\begin{center}
		{\Large\bfseries{\@title}\par}
		\vskip 1em
		{\normalsize
			\lineskip .5em
			\begin{tabular}[t]{c}
				\@author
			\end{tabular}\par}
		\vskip 1.5em
	\end{center}
}
\makeatother

\usepackage{color}
\usepackage{xcolor}
\usepackage[normalem]{ulem} 
\usepackage{soul}

\usepackage{graphicx} % already included
\usepackage{xstring}  % for string manipulation

\usepackage{xstring}
\usepackage{graphicx} % for \scalebox

\usepackage{marvosym}
\begin{document}
	\begin{center}{\bf \large
			{Certain Congruences for Simultaneously $s$-Regular and $t$-Distinct Partition \\[0.15cm] }  }\vspace{0.2cm}
	\end{center}
	\begin{center}
		\footnotemark[1]
		\bf Anjelin Mariya Johnson, \footnotemark[2]
		\bf S. N. Fathima  %\footnotetext{Email: dr.fathima.sn@pondiuni.ac.in}

	\end{center}
	
	%%%%%%%%%%%%%%%%%%%%%%%%%%%%%%%%%%%%%%%%%%%%%%%%%%%%%%%%%%%%%%%%%%%
	\vspace{1 cm}
	%\begin{abstract}
	\begin{center}
		\begin{minipage}{0.85\textwidth}  % abstract uses only 80% of the page width
			
			\noindent{\bf Abstract:} 
			Let $RD^{\ell,t}(n)$ denote the number of simultaneously $l$-regular and  $t$-distinct partitions. In this paper, we prove some infinte families of congruences modulo powers of 5 for $RD^{\ell,t}(n)$. For example, 
			for each $n\geq0$, $0\leq\alpha\leq{s-k-1}$ and $1\leq k <s$,
			\begin{align*}
				\sum_{n=0}^{\infty}RD_{(5^{2k-1},5^{2s-1})}\big(5^{2k+2\alpha-1}n+\frac{4\cdot5^{2k+2\alpha}-5^{2k-1}+1}{24}\big)q^n 
				\equiv 0 \pmod{5^{2k+\alpha-1}}
			\end{align*}
			\vspace{.125cm}
			
			\vspace{.05cm}
			
			%\end{abstract}
		\end{minipage}
	\end{center}
	%%%%%%%%%%%%%%%%%%%%%%%%%%%%%%%%%%%%%%%%%%%%%%%%%%%%%%%%%%%%%%%%%%%
	
	\section {Introduction}
A partition $\lambda=(\lambda_1,\lambda_2,\ldots,\lambda_k)$ of a non-negative integer $n$ is a non-increasing sequence of positive integers summing to $n$. Each $\lambda_i$ is called a part. Let the number of partitions of $n$ be denoted by $p(n)$. By convention, $p(0):=1$. The generating function for $p(n)$ is given by
\begin{align*}
	\sum_{n=0}^{\infty}p(n)q^n=\prod_{n=1}^{\infty}\frac{1}{1-q^n}=\frac{1}{(q;q)_\infty}
\end{align*}
where 
\begin{align*}
	(a;q)_\infty:=\prod_{i=0}^{\infty}(1-aq^i),\qquad|q|<1.
\end{align*}
Throughout the paper, we use the notation $E_k:=(q^k;q^k)_\infty$.
	%%%%%%%%%%%%%%%%%%%%%%%%%%%%%%%%%%%%%%%%%%%%%%%%%%%%%%%%%%%%%%%%%%%%%%
	The arithmetic properties of restricted partition families have been a central focus in number theory. A prominent family in restricted partition is simultaneously $l$-regular and $t$-distinct partition function, $RD^{\ell,t}(n)$, introduced by Keith \cite{keith}, where none of the parts are divisible by $l$ and parts appear fewer than $t$ times. The generating function for $RD^{\ell,t}$ is given by
	\begin{align}\label{Gf}
		\sum_{n=0}^{\infty} RD^{\ell,t}(n)q^n=
		\frac{E_{\ell}E_t}{E_1E_{\ell t}}.
	\end{align}
	For example, four partitions of $RD^{3,3}(6)$ are: $(5,1),\;(4,2),\;(4,1,1),\;(2,2,1,1).$\\
Congruences for $RD_{l,t}(n)$ have so far been obtained mainly for small moduli. Buragohain and Saikia \cite{saikia} studied several pairs $(l,t)$ including $(5\alpha,5\beta)$ and obtained the following result
\begin{theorem}[Buragohain and Saikia {\cite[Theorem 5]{saikia}}]
	For any integers $\alpha > 0$, $\beta > 0$ and prime $p$, we have
	\begin{align*}
		RD_{5\alpha,5\beta}(5n+4) &\equiv 0 \pmod{5},\\
		RD_{p,p}(n) \equiv& \delta_{(p-1)^{2}}(n) \pmod{p}.
	\end{align*}
	where $\delta_{m}(n)$ denotes the number of partitions of a positive
	integer $n$ with each part having $m$ different colors.
\end{theorem}
 For more information on this partition function $RD^{\ell,t}(n)$, we refer the readers to \cite{keith,keith2,Ahmia 2}.\\
In this paper, we revisit the partition function $RD^{\ell,t}(n)$ and study the pair $(5^{2k-1},5^{2s-1})$. We prove infinite families of congruences modulo powers of $5$, where the exponent of the modulus is given explicitly in terms of $k$, $s$ and $\alpha$.  
Our main results are stated below.
	
	%%%%%%%%%%%%%%%%%%%%%%%%%%%%%%%%%%%%%%%%%%%%%%%%%%%%%%%%%%%%%%%%%%%%%%%%%%%%%%%%%%%%%%%%%
	\noindent The theorem below provides infinte families of congruences modulo powers of 5.
	\begin{theorem}\label{T1}
		For each $n\geq0$, $0\leq\alpha<{s-k-1}$ and $1\leq k \leq s$,
		\begin{align}
			\sum_{n=0}^{\infty}RD_{(5^{2k-1},5^{2s-1})}\big(5^{2k+2\alpha-1}n+\frac{4\cdot5^{2k+2\alpha}-5^{2k-1}+1}{24}\big)q^n 
			\equiv 0 \pmod{5^{2k+\alpha-1}}\label{t11}\\
				\sum_{n=0}^{\infty}RD_{(5^{2k-1},5^{2s-1})}\big(5^{2k+2\alpha}n+\frac{(20+24r)\cdot5^{2k+2\alpha-1}-5^{2k-1}+1}{24}\big)q^n 
			\equiv 0 \pmod{5^{2k+\alpha}},\label{t12} 
		\end{align}
		\text{where $r\in \{1,2,3,4\}$}.
	\end{theorem}
	\begin{theorem}\label{t2}
		For each $n\geq0$, $0\leq\alpha\leq{k-1}$ and $1\leq k < s$,
		\begin{align}
			\sum_{n=0}^{\infty}RD_{(5^{2k-1},5^{2s-1})}\Big(5^{2s+2\alpha-1}n+\frac{21\cdot5^{2s+2\alpha-1}-5^{2s-1}-5^{2k-1}+1}{24}\Big)q^n\nonumber\\
			\equiv0\pmod{5^{k+s+\alpha-1}}\label{t21},\\
		\sum_{n=0}^{\infty}RD_{(5^{2k-1},5^{2s-1})}\Big(5^{2s+2\alpha}n+\frac{9\cdot5^{2s+2\alpha}-5^{2s-1}-5^{2k-1}+1}{24}\Big)q^n\nonumber\\
		\equiv0\pmod{5^{k+s+\alpha-1}}\label{t22},\\
			\sum_{n=0}^{\infty}RD_{(5^{2k-1},5^{2s-1})}\Big(5^{2s+2\alpha}n+\frac{93\cdot5^{2s+2\alpha-1}-5^{2s-1}-5^{2k-1}+1}{24}\Big)q^n\nonumber\\
		\equiv0\pmod{5^{k+s+\alpha}}\label{t23},\\
			\sum_{n=0}^{\infty}RD_{(5^{2k-1},5^{2s-1})}\Big(5^{2s+2\alpha}n+\frac{117\cdot5^{2s+2\alpha-1}-5^{2s-1}-5^{2k-1}+1}{24}\Big)q^n\nonumber\\
		\equiv0\pmod{5^{k+s+\alpha}}\label{t24}.
		\end{align}
	\end{theorem}
	\begin{theorem}\label{t3}
		For each $n\geq0$, $\alpha\geq{0}$ and $1\leq k < s$,
		\begin{align}
				\sum_{n=0}^{\infty}RD_{(5^{2k-1},5^{2s-1})}\Big(5^{2s+2k+2\alpha-3}n+\frac{16\cdot5^{2s+2k+2\alpha-3}+5^{2s+2k-2}-5^{2s-1}-5^{2k-1}+1}{24}\Big)q^n\nonumber\\
		\equiv0\pmod{5^{2k+s+\alpha-2}}\label{t31},\\
			\sum_{n=0}^{\infty}RD_{(5^{2k-1},5^{2s-1})}\Big(5^{2s+2k+2\alpha-2}n+\frac{8\cdot5^{2s+2k+2\alpha-2}+5^{2s+2k-2}-5^{2s-1}-5^{2k-1}+1}{24}\Big)q^n\nonumber\\
			\equiv0\pmod{5^{2k+s+\alpha-2}}\label{t32},\\
			\sum_{n=0}^{\infty}RD_{(5^{2k-1},5^{2s-1})}\Big(5^{2s+2k+2\alpha-3}n+\frac{(16+24r)\cdot5^{2s+2k+2\alpha-3}+5^{2s+2k-2}-5^{2s-1}-5^{2k-1}+1}{24}\Big)q^n\nonumber\\
			\equiv0\pmod{5^{2k+s+\alpha-1}}\label{t33},\qquad r\in \{3,4\}\\
				\sum_{n=0}^{\infty}RD_{(5^{2k-1},5^{2s-1})}\Big(5^{2s+2k+2\alpha-2}n+\frac{(8+24r)\cdot5^{2s+2k+2\alpha-2}+5^{2s+2k-2}-5^{2s-1}-5^{2k-1}+1}{24}\Big)q^n\nonumber\\
			\equiv0\pmod{5^{2k+s+\alpha-1}},\qquad r\in \{2,3,4\}\label{t34}.
			\end{align}
	\end{theorem}
	%%%%%%%%%%%%%%%%%%%%%%%%%%%%%%%%%%%%%%%%%%%%%%%%%%%%%%%%%%%%%%%%%%%%%%%%%
	\noindent	We organise the paper as follows: Section 2 recalls the preliminary results that the later arguments rely on. In Section 3 we obtain explicit generating functions for $RD(n)$ along arithmetic progressions. Section 4 estimates the 5-adic valuations of these coefficients and prove Theorems 1.2–1.4 with the help of $H$ operator due to Hirschhorn and Hunt.
	\section{Preliminaries} 
	In this section, we list certain definitions and lemmas that plays a vital role in proving our results.\\
	Using Binomial theorem, we obtain a congruence result which is used frequently in the sequel without explicitly mentioning.
	\begin{align}
		E_m^{p^k} \equiv E_{mp}^{p^{k-1}} \pmod{p^k}.\label{binomial}
	\end{align}
	%%%%%%%%%%%%%%%%%%%%%%%%%%%%%%%%%%%%%%%%%%%%%%%%%%%%%%%%%%%%%%%%%
	We next recall the famous Rogers Ramanujan continued fraction $R(q)$\cite{power}:
	\begin{align*}
		R(q):=\frac{(q;q^5)_\infty(q^4;q^5)_\infty}{(q^2;q^5)_\infty(q^3;q^5)_\infty}.
	\end{align*}
	which satisfies the following 5 dissection for $E_1$ involving $R(q)$
	\begin{align}
		E_1=E_{25}\left(\frac{1}{R(q^5)}-q-q^2R(q^5)\right).\label{r(q)}
	\end{align}
	Thanks to (\ref{r(q)}) we now define
	\begin{align*}
		\zeta = \frac{E_1}{qE_{25}}, \qquad T = \frac{E_5^6}{q^5 E_{25}^6}.
	\end{align*}
	
	Following Hirschhorn \cite{power, Hunt}, we introduce a ``huffing'' operator $H$, which operates on a series by picking out those terms of the form $q^{5n}$, and huffing the rest away. That is,
	\[
	H\left( \sum_{n=0}^{\infty} a_n q^n \right) = \sum_{n=0}^{\infty} a_{5n} q^{5n}.
	\]
	As in Hirschhorn \cite{power}, we define an infinite matrix $\{m(i,j)\}_{i,j\geq 1}$ by
	\[
	\begin{pmatrix}
		5 & 0 & 0 & 0 & 0 & 0 & \cdots\\
		2\times 5 & 5^3 & 0 & 0 & 0 & 0 & \cdots\\
		9 & 3\times 5^3 & 5^5 & 0 & 0 & 0 & \cdots\\
		4 & 22\times 5^2 & 4\times 5^5 & 5^7 & 0 & 0 & \cdots\\
		1 & 4\times 5^3 & 8\times 5^5 & 5^8 & 5^9 & 0 & \cdots\\
		\vdots & \vdots & \vdots & \vdots & \vdots & \vdots & \ddots
	\end{pmatrix}
	\]
	and for $i\geq 6$, $m(i,1)=0$, and for $j\geq 2$,
	\begin{align}
		m(i,j) &= 25m(i-1,j-1) + 25m(i-2,j-1) + 15m(i-3,j-1) \notag\\
		&\qquad + 5m(i-4,j-1) + m(i-5,j-1). \label{eq:2.4}
	\end{align}
	By induction, it follows immediately that
	
	\begin{lemma}[\cite{power}] \label{lem:2.1}
		We have
		\begin{itemize}
			\item $m(i,j)=0$ for $j>i$.
			\item $m(i,j)=0$ for $i>5j$.
		\end{itemize}
	\end{lemma}
	
	The following lemmas play a vital role to prove our main results.
	
	\begin{lemma}[\cite{power}, Eq.~(6.4.9)] \label{lem:2.2}
		For $j\geq 1$, we have
		\[
		H\left(\frac{1}{\zeta^i}\right) = \sum_{j=1}^{\infty} \frac{m(i,j)}{T^j} = \sum_{j=1}^{i} \frac{m(i,j)}{T^j}.
		\]
	\end{lemma}
	
	%%%%%%%%%%%%%%%%%%%%%%%%%%%%%%%%%%%%%%%%%%%%%%%%%%%%%%%%%%%%%%%%
	%%%%%%%%%%%%%%%%%%%%%%%%%%%%%%%%%%%%%%%%%%%%%%%%%%%%%%%%%%%%%%%%%%%%%%%%%%%%%%%%%%%%%%%%%
	\begin{lemma}
		For integers $i,j \geq 1$, we have
	\begin{align}
		H\!\left(q^{i-1}\frac{E_{5}^{6i-1}}{E_{1}^{6i-1}}\right)
		&=\sum_{j=1}^{5i} a_{i,j}\,q^{5j-5}\frac{E_{25}^{6j-5}}{E_{5}^{6j-5}},\label{L1}\\
		H\!\left(q^{i-5}\Big(\frac{E_{5}}{E_{1}}\Big)^{6i-5}\right)
		&=\sum_{j=1}^{5i-4} b_{i,j}\,q^{5j-5}\Big(\frac{E_{25}}{E_{5}}\Big)^{6j-1},\label{L2}\\
		H\!\left(q^{i-2}\frac{E_5^{6i-1}}{E_1^{6i-2}}\right)
		&=\sum_{j=1}^{5i-1} c_{i,j}\,q^{5j-5}\frac{E_{25}^{6j-4}}{E_{5}^{6j-5}},\label{L3}\\
		H\!\left(q^{i-5}\frac{E_5^{6i-4}}{E_1^{6i-5}}\right)
		&=\sum_{j=1}^{5i-4} b_{i,j}\,q^{5j-5}\frac{E_{25}^{6j-1}}{E_{5}^{6j-2}},\label{L4}\\
		H\!\left(q^{i-2}\Big(\frac{E_5}{E_1}\Big)^{6i-2}\right)
		&=\sum_{j=1}^{5i-1} c_{i,j}\,q^{5j-5}\Big(\frac{E_{25}}{E_{5}}\Big)^{6j-4},\label{L5}\\
		H\!\left(q^{i-4}\Big(\frac{E_5}{E_1}\Big)^{6i-4}\right)
		&=\sum_{j=1}^{5i-3} d_{i,j}\,q^{5j-5}\Big(\frac{E_{25}}{E_{5}}\Big)^{6j-2}.\label{L6}
	\end{align}
	\end{lemma}
	\begin{proof}
		Note for $1\leq j \leq i-1$, $m_{6i-1,j}=0$, therefore, we have
		\noindent
		\begin{align*}\label{a15}
			H\!\left(\frac{1}{\zeta^{6i-1}}\right)
			= \sum_{j= i}^{6i-1} \frac{m_{6i-1,\,j}}{T^{j}}
			= \sum_{j= 1}^{5i} \frac{m_{6i-1,\,i+j-1}}{T^{i+j-1}},
		\end{align*}
		which can be rewritten as
		\begin{equation*}
			H\!\left(q^{6i-1}\frac{E_{25}^{6i-1}}{E_{1}^{6i-1}}\right)
			=
			\sum_{j=1}^{5i}
			m_{6i-1,i+j-1}\,
			q^{5i+5j-5}
			\frac{E_{25}^{6i+6j-6}}{E_{5}^{6i+6j-6}},
		\end{equation*}
		that gives
		\begin{equation*}
			H\!\left(q^{i-1}\frac{E_{5}^{6i-1}}{E_{1}^{6i-1}}\right)=
		\sum_{j=1}^{5i} a_{i,j}\,
		q^{5j-5}\frac{E_{25}^{6j-5}}{E_{5}^{6j-5}},\\
		\end{equation*}
		where
		\[
		a_{i,j}=m_{6i-1,i+j-1},
		\]
		this completes the proof of (\ref{L1}).\\
		Similarly, for $1\leq j \leq i-1$, $m_{6i-5,j}=0$, then
		\begin{align*}
			H\!\left(\frac{1}{\zeta^{6i-5}}\right)
			= \sum_{j= i}^{6i-5} \frac{m_{6i-5,\,j}}{T^{j}}
			= \sum_{j= 1}^{5i-4} \frac{m_{6i-5,\,i+j-1}}{T^{i+j-1}},
		\end{align*}
		i.e,
		\begin{equation*}
			H\!\left(q^{6i-5}\frac{E_{25}^{6i-5}}{E_{1}^{6i-5}}\right)
			=
			\sum_{j=1}^{5i-4}
			m_{6i-5,i+j-1}\,
			q^{5i+5j-5}
			\frac{E_{25}^{6i+6j-6}}{E_{5}^{6i+6j-6}},
		\end{equation*}
		that gives
		\begin{equation*}
			H\!\left(q^{i-5}{(\frac{E_{5}}{E_{1}})}^{6i-5}\right)=
			\sum_{j=1}^{5i-4} b_{i,j}\,
			q^{5j-5}{(\frac{E_{25}}{E_{5}})}^{6j-1}. 
		\end{equation*}
		where
		\[
		b_{i,j}=m_{6i-5,i+j-1},
		\]
		this completes the proof of (\ref{L2}). Similar proof follows for (\ref{L3})-(\ref{L6})
	\end{proof}
	\section{Generating Functions}
In this section, we esatblish the generating functions for $RD^{\ell,t}(n)$ along certain arithmetic progressions. These generating functions are the main tool for the congruences in the next section.
	\begin{theorem}\label{gf1}
		For each $n\geq0$, $0\leq\alpha\leq{s-k-1}$ and $1\leq k < s$,
		\begin{align}
			\sum_{n=0}^{\infty}RD_{(5^{2k-1},5^{2s-1})}\big(5^{2k+2\alpha-1}n+\frac{4\cdot5^{2k+2\alpha}-5^{2k-1}+1}{24}\big)q^n \nonumber\\=\frac{{E_{5^{2(s-k)-2\alpha}}}}{{E_{5^{2s-2\alpha-1}}}}\sum_{i=1}^{\infty}y_{2\alpha+1,i}^{2k-1}q^{i-1}\frac{E_5^{6i-1}}{E_1^{6i-1}},\label{gf11}\\
			\sum_{n=0}^{\infty}RD_{(5^{2k-1},5^{2s-1})}\big(5^{2k+2\alpha}n+\frac{4\cdot5^{2k+2\alpha}-5^{2k-1}+1}{24}\big)q^n \nonumber\\=\frac{{E_{5^{2(s-k)-2\alpha-1}}}}{{E_{5^{2s-2\alpha-2}}}}\sum_{i=1}^{\infty}y_{2\alpha+2,i}^{2k-1}q^{i-1}\frac{E_5^{6i-5}}{E_1^{6i-5}}\label{gf12}.
		\end{align}
		where the coefficient entries are defined as follows:
		\begin{align*}
			Y_1^{2k-1}=&X_{2k-1}
		\end{align*}
		and for any $1\leq\alpha\leq{s-k-1}$,
		\begin{align}\label{Ydef}
			Y_{\alpha+1}^{2k-1}=
			\begin{cases}
				Y_\alpha^{2k-1}\cdot A, & \text{if $\alpha$ is odd}\\
				Y_\alpha^{2k-1}\cdot B, & \text{if $\alpha$ is even}.
			\end{cases}		
		\end{align}
			where,
		\begin{align*}
			A:=&(a_{i,j})=m_{6i-1,i+j-1}\\
			B:=&(b_{i,j})=m_{6i-5,i+j-1}.
		\end{align*}
		The coefficient vectors $X_1=(5,0,0,0,\ldots)$,
		and for any $t\geq1,$ we have\\
		\begin{align*}
			X_{t+1}=\begin{cases}
				 X_t \cdot P,& \text{if $t$ is odd}\\
				X_t \cdot Q,& \text{if $t$ is even}
			\end{cases}
		\end{align*}
		where,
		\begin{align*}
			P:=&(p_{i,j})=m_{6i,i+j}\\
			Q:=&(q_{i,j})=m_{6i+1,i+j}.
		\end{align*}
	\end{theorem}
	%%%%%%%%%%%%%%%%%%%%%%%%%%%%%%%%%%%%%%%%%%%%%%%%%%%%%%%%%%

	%%%%%%%%%%%%%%%%%%%%%%%%%%%%%%%%%%%%%%%%%%%%%%%%%%%%%%%%%%%%
	\begin{proof}
		It follows from (\ref{Gf}) that
		\begin{align}\label{3.1.1}
			\sum_{n=0}^{\infty} RD_{(5^{2k-1},5^{2s-1})}(n)\,q^n=\frac{{E_{5^{2k-1}}\cdot E_{5^{2s-1}}}}{{E_1}{E_{5^{2(k+s-1)}}}}
			=\frac{{E_{5^{2k-1}}\cdot E_{5^{2s-1}}}}{{E_{5^{2(k+s-1)}}}} \sum_{n=0}^{\infty} p(n)\,q^n,
		\end{align}
		where $p(n)$ is defined by
		\[
		\sum_{n=0}^{\infty} p(n)\,q^n = \frac{1}{{E_1}}.
		\]
		Extracting the terms involving $q^{5^{2k-1}n+(19\times 5^{2k-1}+1)/24}$ in (\ref{3.1.1}) and upon simplification, gives
		\begin{align}
			\sum_{n=0}^{\infty} RD_{(5^{2k-1},5^{2s-1})}\!\left(5^{2k-1}n+\frac{19\times 5^{2k-1}+1}{24}\right)q^n\nonumber\\
			= \frac{E_1{E_5}^{2(s-k)}}{{E_{5^{2s-1}}}^2} \sum_{n=0}^{\infty} p\!\left(5^{2k-1}n+\frac{19\times 5^{2k-1}+1}{24}\right)q^n. \label{after}
		\end{align}
		Thanks to Hirschhorn et al. \cite[Theorem (1.4)]{Hunt}, we have
		\begin{align}
			\sum_{n=0}^{\infty} p\!\left(5^{2k-1}n+\frac{19\times 5^{2k-1}+1}{24}\right)q^n
			= \sum_{j=1}^{\infty} x_{2k-1,j}\, q^{j-1}\,\frac{E_5^{6j-1}}{E_1^{6j}}.\label{tang}
		\end{align}
		Employing (\ref{tang} ) in (\ref{after}) yields
		\[
		\sum_{n=0}^{\infty} RD_{(5^{2k-1},5^{2s-1})}\!\left(5^{2k-1}n+\frac{19\times 5^{2k-1}+1}{24}\right)q^n
		= \frac{{E_5}^{2(s-k)}}{{E_{5^{2s-1}}}^2} \sum_{j=1}^{\infty} x_{2k-1,j}\, q^{j-1}\,\frac{E_5^{6j-1}}{E_1^{6j-1}}.
		\]
		This implies $\alpha=0$ case of Theorem \ref{gf1} holds.
		
		Now for some $\alpha\ge 0$, we assume that Theorem \ref{gf1} holds, specifically
		\begin{align*}
			\sum_{n=0}^{\infty}RD_{(5^{2k-1},5^{2s-1})}\big(5^{2k+2\alpha-1}n+\frac{4\cdot5^{2k+2\alpha}-5^{2k-1}+1}{24}\big)q^n \nonumber\\=\frac{{E_{5^{2(s-k)-2\alpha}}}}{{E_{5^{2s-2\alpha-1}}}}\sum_{i=1}^{\infty}y_{2\alpha+1,i}^{2k-1}q^{i-1}\frac{E_5^{6i-1}}{E_1^{6i-1}}.\\
		\end{align*}
		Applying the $H$-operator to both sides of the above equation, we obtain
		\begin{align*}
			\sum_{n=0}^{\infty}RD_{(5^{2k-1},5^{2s-1})}\big(5^{2k+2\alpha-1}(5n)&+\frac{4\cdot5^{2k+2\alpha}-5^{2k-1}+1}{24}\big)q^{5n} \nonumber\\=&\frac{{E_{5^{2(s-k)-2\alpha}}}}{{E_{5^{2s-2\alpha-1}}}}\sum_{i=1}^{\infty}y_{2\alpha+1,i}^{2k-1}H\big(q^{i-1}\frac{E_5^{6i-1}}{E_1^{6i-1}}\big)\\
		=&\frac{{E_{5^{2(s-k)-2\alpha}}}}{{E_{5^{2s-2\alpha-1}}}}\sum_{i=1}^{\infty}y_{2\alpha+1,i}^{2k-1}\sum_{j=1}^{5i} a_{i,j}\,
		q^{5j-5}\frac{E_{25}^{6j-5}}{E_{5}^{6j-5}}\\
		=&\frac{{E_{5^{2(s-k)-2\alpha}}}}{{E_{5^{2s-2\alpha-1}}}}\sum_{j=1}^{\infty}\big(\sum_{i=1}^{\infty}y_{2\alpha+1,i}^{2k-1}a_{i,j}\big)q^{5j-5}\frac{E_{25}^{6j-5}}{E_{5}^{6j-5}}\\
		=&\frac{{E_{5^{2(s-k)-2\alpha}}}}{{E_{5^{2s-2\alpha-1}}}}\sum_{j=1}^{\infty}y_{2\alpha+2,i}^{2k-1}q^{5j-5}\frac{E_{25}^{6j-5}}{E_5^{6j-5}}.
		\end{align*}
			Finally, on	replacing $q^5$ by $q$ we have
			\begin{align*}
				\sum_{n=0}^{\infty}RD_{(5^{2k-1},5^{2s-1})}\big(5^{2k+2\alpha}n+\frac{4\cdot5^{2k+2\alpha}-5^{2k-1}+1}{24}\big)q^n \nonumber\\=\frac{{E_{5^{2(s-k)-2\alpha-1}}}}{{E_{5^{2s-2\alpha-2}}}}\sum_{i=1}^{\infty}y_{2\alpha+2,i}^{2k-1}q^{i-1}\frac{E_5^{6i-5}}{E_1^{6i-5}}.
			\end{align*}
			which is (\ref{gf12}).\\
				Assume that (\ref{gf12}) is true for some $\alpha\geq0$.
				\begin{align*}
					\sum_{n=0}^{\infty}RD_{(5^{2k-1},5^{2s-1})}\big(5^{2k+2\alpha}n+\frac{4\cdot5^{2k+2\alpha}-5^{2k-1}+1}{24}\big)q^{n-4} \nonumber\\=\frac{{E_{5^{2(s-k)-2\alpha-1}}}}{{E_{5^{2s-2\alpha-2}}}}\sum_{i=1}^{\infty}y_{2\alpha+2,i}^{2k-1}q^{i-5}\frac{E_5^{6i-5}}{E_1^{6i-5}}.
				\end{align*}
					Applying $H$-operator gives
					\begin{align*}
						\sum_{n=0}^{\infty}RD_{(5^{2k-1},5^{2s-1})}\big(5^{2k+2\alpha}(5n+&4)+\frac{4\cdot5^{2k+2\alpha}-5^{2k-1}+1}{24}\big)q^{5n} \nonumber\\=&\frac{{E_{5^{2(s-k)-2\alpha-1}}}}{{E_{5^{2s-2\alpha-2}}}}\sum_{i=1}^{\infty}y_{2\alpha+2,i}^{2k-1}H\big(q^{i-5}\frac{E_5^{6i-5}}{E_1^{6i-5}}\big)\\
						=&\frac{{E_{5^{2(s-k)-2\alpha-1}}}}{{E_{5^{2s-2\alpha-2}}}}\sum_{i=1}^{\infty}y_{2\alpha+2,i}^{2k-1}\sum_{j=1}^{5i-4} b_{i,j}q^{5j-5}{(\frac{E_{25}}{E_{5}})}^{6j-1}\\
						=&\frac{{E_{5^{2(s-k)-2\alpha-1}}}}{{E_{5^{2s-2\alpha-2}}}}\sum_{j=1}^{\infty}\big(\sum_{i=1}^{\infty}y_{2\alpha+2,i}^{2k-1} b_{i,j}\big)q^{5j-5}{(\frac{E_{25}}{E_{5}})}^{6j-1}\\
						=&\frac{{E_{5^{2(s-k)-2\alpha-1}}}}{{E_{5^{2s-2\alpha-2}}}}\sum_{j=1}^{\infty}y_{2\alpha+3,j}^{2k-1}q^{5j-5}\frac{E_{25}^{6j-1}}{E_5^{6j-1}}.
					\end{align*}
					Replacing $q^5$ by $q$ 
					\begin{align*}
						\sum_{n=0}^{\infty}RD_{(5^{2k-1},5^{2s-1})}\big(5^{2k+2\alpha+1}n+\frac{4\cdot5^{2k+2\alpha+2}-5^{2k-1}+1}{24}\big)q^n \nonumber\\=\frac{{E_{5^{2(s-k)-2\alpha-2}}}}{{E_{5^{2s-2\alpha-3}}}}\sum_{i=1}^{\infty}y_{2\alpha+3,i}^{2k-1}q^{i-1}\frac{E_5^{6i-1}}{E_1^{6i-1}}.
					\end{align*}
		This implies Theorem \ref{gf1} is true for $\alpha+1$, hence by induction, we complete the proof of Theorem \ref{gf1}.
	\end{proof}
		\begin{theorem}\label{gf2}
		For each $n\geq0$, $0\leq\alpha\leq{k-1}$ and $1\leq k < s$,
		\begin{align}
			\sum_{n=0}^{\infty}RD_{(5^{2k-1},5^{2s-1})}\Big(5^{2s+2\alpha-1}n+\frac{21\cdot5^{2s+2\alpha-1}-5^{2s-1}-5^{2k-1}+1}{24}\Big)q^n\nonumber\\
			=\frac{1}{E_{5^{2k-2\alpha-1}}}\sum_{i=1}^{\infty}z_{2\alpha+1,i}\,q^{i-1}\frac{E_5^{6i-1}}{E_1^{6i-2}},\label{gf21}\\
			\sum_{n=0}^{\infty}RD_{(5^{2k-1},5^{2s-1})}\Big(5^{2s+2\alpha}n+\frac{9\cdot5^{2s+2\alpha}-5^{2s-1}-5^{2k-1}+1}{24}\Big)q^n\nonumber\\
			=\frac{1}{E_{5^{2k-2\alpha-2}}}\sum_{i=1}^{\infty}z_{2\alpha+2,i}\,q^{i-1}\frac{E_5^{6i-4}}{E_1^{6i-5}}.\label{gf22}
		\end{align}
	where the coefficient vectors are defined as follows:
		\[
		Z_1=Y_{2s-2k}^{2k-1}B,
		\]
		and for $1\leq\beta\leq 2k-1$,
		\begin{align}\label{zdef}
			Z_{\beta+1}=\begin{cases}
			Z_\beta\cdot C, & \text{if $\beta$ is odd},\\[2pt]
			Z_\beta\cdot B, & \text{if $\beta$ is even},
		\end{cases}
		\end{align}
		where, \begin{align*}
			C:=(c_{i,j})=m_{6i-2,i+j-1},\\
			 B:=(b_{i,j})=m_{6i-5,i+j-1}.
		\end{align*}
	\end{theorem}
	\begin{proof}
		For case, $\alpha=s-k-1$ in Theorem \ref{gf1} gives
		\begin{align*}
			\sum_{n=0}^{\infty}RD_{(5^{2k-1},5^{2s-1})}\Big(5^{2s-2}n+\frac{4\cdot5^{2s-2}-5^{2k-1}+1}{24}\Big)q^{n-4}\\
			=\frac{E_5}{E_{5^{2k}}}\sum_{i=1}^{\infty}y_{2s-2k,i}^{2k-1}\,q^{i-5}\frac{E_5^{6i-5}}{E_1^{6i-5}}.
		\end{align*}
		Applying $H$-operator in above identity, we have
		\begin{align*}
			\sum_{n=0}^{\infty}RD_{(5^{2k-1},5^{2s-1})}\Big(5^{2s-2}(5&n+4)+\frac{4\cdot5^{2s-2}-5^{2k-1}+1}{24}\Big)q^{5n}\\
			=&\frac{E_5}{E_{5^{2k}}}\sum_{i=1}^{\infty}y_{2s-2k,i}^{2k-1}H\Big(q^{i-5}\frac{E_5^{6i-5}}{E_1^{6i-5}}\Big)\\
			=&\frac{E_5}{E_{5^{2k}}}\sum_{i=1}^{\infty}y_{2s-2k,i}^{2k-1}\sum_{j=1}^{5i-4} b_{i,j}\,
			q^{5j-5}\Big(\frac{E_{25}}{E_{5}}\Big)^{6j-1}\\
			=&\frac{E_5}{E_{5^{2k}}}\sum_{j=1}^{\infty}\Big(\sum_{i=1}^{\infty}y_{2s-2k,i}^{2k-1}b_{i,j}\Big)q^{5j-5}\Big(\frac{E_{25}}{E_{5}}\Big)^{6j-1}\\
			=&\frac{E_5}{E_{5^{2k}}}\sum_{j=1}^{\infty}z_{1,j}\,q^{5j-5}\Big(\frac{E_{25}}{E_{5}}\Big)^{6j-1}.
		\end{align*}
		Finally, on replacing $q^5$ by $q$ we have
		\begin{align*}
			\sum_{n=0}^{\infty}RD_{(5^{2k-1},5^{2s-1})}\Big(5^{2s-1}n+\frac{21\cdot5^{2s-1}-5^{2s-1}-5^{2k-1}+1}{24}\Big)q^{n}\\
			=\frac{1}{E_{5^{2k-1}}}\sum_{j=1}^{\infty}z_{1,j}\,q^{j-1}\frac{E_5^{6j-1}}{E_1^{6j-2}}
		\end{align*}
		which is the $\alpha=0$ case of (\ref{gf21}).\\
		Now, suppose (\ref{gf21}) is true for some $0\leq\alpha\leq k-1$, we have
		\begin{align*}
			\sum_{n=0}^{\infty}RD_{(5^{2k-1},5^{2s-1})}\Big(5^{2s+2\alpha-1}n+\frac{21\cdot5^{2s+2\alpha-1}-5^{2s-1}-5^{2k-1}+1}{24}\Big)q^{n-1}\\
			=\frac{1}{E_{5^{2k-2\alpha-1}}}\sum_{i=1}^{\infty}z_{2\alpha+1,i}\,q^{i-2}\frac{E_5^{6i-1}}{E_1^{6i-2}}.
		\end{align*}
		Applying $H$-operator to both sides, we have
		\begin{align*}
			\sum_{n=0}^{\infty}RD_{(5^{2k-1},5^{2s-1})}\Big(5^{2s+2\alpha-1}(5n+1)+&\frac{21\cdot5^{2s+2\alpha-1}-5^{2s-1}-5^{2k-1}+1}{24}\Big)q^{5n}\\
			=&\frac{1}{E_{5^{2k-2\alpha-1}}}\sum_{i=1}^{\infty}z_{2\alpha+1,i}\,H\Big(q^{i-2}\frac{E_5^{6i-1}}{E_1^{6i-2}}\Big)\\
			=&\frac{1}{E_{5^{2k-2\alpha-1}}}\sum_{i=1}^{\infty}z_{2\alpha+1,i}\sum_{j=1}^{5i-1} c_{i,j}\,
			q^{5j-5}\frac{E_{25}^{6j-4}}{E_{5}^{6j-5}}\\
			=&\frac{1}{E_{5^{2k-2\alpha-1}}}\sum_{j=1}^{\infty}\Big(\sum_{i=1}^{\infty}z_{2\alpha+1,i}\,c_{i,j}\Big)q^{5j-5}\frac{E_{25}^{6j-4}}{E_{5}^{6j-5}}\\
			=&\frac{1}{E_{5^{2k-2\alpha-1}}}\sum_{j=1}^{\infty}z_{2\alpha+2,j}\,q^{5j-5}\frac{E_{25}^{6j-4}}{E_{5}^{6j-5}}.
		\end{align*}
		Replacing $q^5$ by $q$ gives
		\begin{align*}
			\sum_{n=0}^{\infty}RD_{(5^{2k-1},5^{2s-1})}\Big(5^{2s+2\alpha}n+\frac{9\cdot5^{2s+2\alpha}-5^{2s-1}-5^{2k-1}+1}{24}\Big)q^n\\
			=\frac{1}{E_{5^{2k-2\alpha-2}}}\sum_{i=1}^{\infty}z_{2\alpha+2,i}\,q^{i-1}\frac{E_5^{6i-4}}{E_1^{6i-5}}
		\end{align*}
		which is (\ref{gf22}).\\
		Suppose that (\ref{gf22}) is true for some $0\leq\alpha\leq k-2$
		\begin{align*}
			\sum_{n=0}^{\infty}RD_{(5^{2k-1},5^{2s-1})}\Big(5^{2s+2\alpha}n+\frac{9\cdot5^{2s+2\alpha}-5^{2s-1}-5^{2k-1}+1}{24}\Big)q^{n-4}\\
			=\frac{1}{E_{5^{2k-2\alpha-2}}}\sum_{i=1}^{\infty}z_{2\alpha+2,i}\,q^{i-5}\frac{E_5^{6i-4}}{E_1^{6i-5}}.
		\end{align*}
		Applying $H$-operator gives
		\begin{align*}
			\sum_{n=0}^{\infty}RD_{(5^{2k-1},5^{2s-1})}\Big(5^{2s+2\alpha}(5n+4)+&\frac{9\cdot5^{2s+2\alpha}-5^{2s-1}-5^{2k-1}+1}{24}\Big)q^{5n}\\
			=&\frac{1}{E_{5^{2k-2\alpha-2}}}\sum_{i=1}^{\infty}z_{2\alpha+2,i}\,H\Big(q^{i-5}\frac{E_5^{6i-4}}{E_1^{6i-5}}\Big)\\
			=&\frac{1}{E_{5^{2k-2\alpha-2}}}\sum_{i=1}^{\infty}z_{2\alpha+2,i}\sum_{j=1}^{5i-4} b_{i,j}\,
			q^{5j-5}\frac{E_{25}^{6j-1}}{E_{5}^{6j-2}}\\
			=&\frac{1}{E_{5^{2k-2\alpha-2}}}\sum_{j=1}^{\infty}\Big(\sum_{i=1}^{\infty}z_{2\alpha+2,i}\,b_{i,j}\Big)q^{5j-5}\frac{E_{25}^{6j-1}}{E_{5}^{6j-2}}\\
			=&\frac{1}{E_{5^{2k-2\alpha-2}}}\sum_{j=1}^{\infty}z_{2\alpha+3,j}\,q^{5j-5}\frac{E_{25}^{6j-1}}{E_{5}^{6j-2}}.
		\end{align*}
		Replacing $q^5$ by $q$
		\begin{align*}
			\sum_{n=0}^{\infty}RD_{(5^{2k-1},5^{2s-1})}\Big(5^{2s+2\alpha+1}n+\frac{21\cdot5^{2s+2\alpha+1}-5^{2s-1}-5^{2k-1}+1}{24}\Big)q^n\\
			=\frac{1}{E_{5^{2k-2\alpha-3}}}\sum_{i=1}^{\infty}z_{2\alpha+3,i}\,q^{i-1}\frac{E_5^{6i-1}}{E_1^{6i-2}}
		\end{align*}
		which is the $\alpha+1$ case of (\ref{gf21}).\\This completes the proof of Theorem \ref{gf2} by induction.
	\end{proof}
\begin{theorem}\label{3.1}
	Let $s>k\geq1$. For each $n\geq0$ and $\alpha\geq0$,
	\begin{align}
		\sum_{n=0}^{\infty}RD_{(5^{2k-1},5^{2s-1})}\Big(5^{2s+2k+2\alpha-3}n+\frac{16\cdot5^{2s+2k+2\alpha-3}+5^{2s+2k-2}-5^{2s-1}-5^{2k-1}+1}{24}\Big)q^n\nonumber\\
		=\sum_{i=1}^{\infty}\tau_{2\alpha+1,i}\,q^{i-1}\Big(\frac{E_5}{E_1}\Big)^{6i-2},\label{gf31}\\
		\sum_{n=0}^{\infty}RD_{(5^{2k-1},5^{2s-1})}\Big(5^{2s+2k+2\alpha-2}n+\frac{8\cdot5^{2s+2k+2\alpha-2}+5^{2s+2k-2}-5^{2s-1}-5^{2k-1}+1}{24}\Big)q^n\nonumber\\
		=\sum_{i=1}^{\infty}\tau_{2\alpha+2,i}\,q^{i-1}\Big(\frac{E_5}{E_1}\Big)^{6i-4}.\label{gf32}
	\end{align}
	where the coefficient vectors are defined as follows:
	\[
	\mathrm{T}_1=Z_{2k-1}=Z_{2k-2}\cdot B,
	\]
	and for $\beta\geq1$,
	\begin{align}
		\mathrm{T}_{\beta+1}=\begin{cases}
			\mathrm{T}_\beta\cdot C, & \text{if $\beta$ is odd},\\[2pt]
			\mathrm{T}_\beta\cdot D, & \text{if $\beta$ is even},\label{tdef}
		\end{cases}
	\end{align}
	where,
	\begin{align*}
		B:=(b_{i,j})=m_{6i-5,i+j-1},\\
		C:=(c_{i,j})=m_{6i-2,i+j-1},\\
		D:=(d_{i,j})=m_{6i-4,i+j-1}.
	\end{align*}
\end{theorem}

\begin{proof}
	For case, $\alpha=k-1$ in (\ref{gf21}) gives
	\begin{align*}
		\sum_{n=0}^{\infty}RD_{(5^{2k-1},5^{2s-1})}\Big(5^{2s+2k-3}n+\frac{21\cdot5^{2s+2k-3}-5^{2s-1}-5^{2k-1}+1}{24}\Big)q^{n}\\
		=\frac{1}{E_{5}}\sum_{i=1}^{\infty}z_{2k-1,i}\,q^{i-1}\frac{E_5^{6i-1}}{E_1^{6i-2}}
		=\sum_{i=1}^{\infty}\tau_{1,i}\,q^{i-1}\Big(\frac{E_5}{E_1}\Big)^{6i-2},
	\end{align*}
	since $21\cdot5^{2s+2k-3}=16\cdot5^{2s+2k-3}+5^{2s+2k-2}$, which is the $\alpha=0$ case of (\ref{gf31}).\\
	Now, suppose (\ref{gf31}) is true for some $\alpha\geq0$, we have
	\begin{align*}
		\sum_{n=0}^{\infty}RD_{(5^{2k-1},5^{2s-1})}\Big(5^{2s+2k+2\alpha-3}n+\frac{16\cdot5^{2s+2k+2\alpha-3}+5^{2s+2k-2}-5^{2s-1}-5^{2k-1}+1}{24}\Big)q^{n-1}\\
		=\sum_{i=1}^{\infty}\tau_{2\alpha+1,i}\,q^{i-2}\Big(\frac{E_5}{E_1}\Big)^{6i-2}.
	\end{align*}
	Applying $H$-operator to both sides, we have
	\begin{align*}
		\sum_{n=0}^{\infty}RD_{(5^{2k-1},5^{2s-1})}\Big(5^{2s+2k+2\alpha-3}(5n+1)+&\frac{16\cdot5^{2s+2k+2\alpha-3}+5^{2s+2k-2}-5^{2s-1}-5^{2k-1}+1}{24}\Big)q^{5n}\\
		=&\sum_{i=1}^{\infty}\tau_{2\alpha+1,i}\,H\Big(q^{i-2}\Big(\frac{E_5}{E_1}\Big)^{6i-2}\Big)\\
		=&\sum_{i=1}^{\infty}\tau_{2\alpha+1,i}\sum_{j=1}^{5i-1} c_{i,j}\,
		q^{5j-5}\Big(\frac{E_{25}}{E_{5}}\Big)^{6j-4}\\
		=&\sum_{j=1}^{\infty}\Big(\sum_{i=1}^{\infty}\tau_{2\alpha+1,i}\,c_{i,j}\Big)q^{5j-5}\Big(\frac{E_{25}}{E_{5}}\Big)^{6j-4}\\
		=&\sum_{j=1}^{\infty}\tau_{2\alpha+2,j}\,q^{5j-5}\Big(\frac{E_{25}}{E_{5}}\Big)^{6j-4}.
	\end{align*}
	Replacing $q^5$ by $q$ gives
	\begin{align*}
		\sum_{n=0}^{\infty}RD_{(5^{2k-1},5^{2s-1})}\Big(5^{2s+2k+2\alpha-2}n+\frac{8\cdot5^{2s+2k+2\alpha-2}+5^{2s+2k-2}-5^{2s-1}-5^{2k-1}+1}{24}\Big)q^n\\
		=\sum_{i=1}^{\infty}\tau_{2\alpha+2,i}\,q^{i-1}\Big(\frac{E_5}{E_1}\Big)^{6i-4}
	\end{align*}
	which is (\ref{gf32}).\\
	Suppose that (\ref{gf32}) is true for some $\alpha\geq0$
	\begin{align*}
		\sum_{n=0}^{\infty}RD_{(5^{2k-1},5^{2s-1})}\Big(5^{2s+2k+2\alpha-2}n+\frac{8\cdot5^{2s+2k+2\alpha-2}+5^{2s+2k-2}-5^{2s-1}-5^{2k-1}+1}{24}\Big)q^{n-3}\\
		=\sum_{i=1}^{\infty}\tau_{2\alpha+2,i}\,q^{i-4}\Big(\frac{E_5}{E_1}\Big)^{6i-4}.
	\end{align*}
	Applying $H$-operator gives
	\begin{align*}
		\sum_{n=0}^{\infty}RD_{(5^{2k-1},5^{2s-1})}\Big(5^{2s+2k+2\alpha-2}(5n+3)+&\frac{8\cdot5^{2s+2k+2\alpha-2}+5^{2s+2k-2}-5^{2s-1}-5^{2k-1}+1}{24}\Big)q^{5n}\\
		=&\sum_{i=1}^{\infty}\tau_{2\alpha+2,i}\,H\Big(q^{i-4}\Big(\frac{E_5}{E_1}\Big)^{6i-4}\Big)\\
		=&\sum_{i=1}^{\infty}\tau_{2\alpha+2,i}\sum_{j=1}^{5i-3} d_{i,j}\,
		q^{5j-5}\Big(\frac{E_{25}}{E_{5}}\Big)^{6j-2}\\
		=&\sum_{j=1}^{\infty}\Big(\sum_{i=1}^{\infty}\tau_{2\alpha+2,i}\,d_{i,j}\Big)q^{5j-5}\Big(\frac{E_{25}}{E_{5}}\Big)^{6j-2}\\
		=&\sum_{j=1}^{\infty}\tau_{2\alpha+3,j}\,q^{5j-5}\Big(\frac{E_{25}}{E_{5}}\Big)^{6j-2}.
	\end{align*}
	Replacing $q^5$ by $q$
	\begin{align*}
		\sum_{n=0}^{\infty}RD_{(5^{2k-1},5^{2s-1})}\Big(5^{2s+2k+2\alpha-1}n+\frac{16\cdot5^{2s+2k+2\alpha-1}+5^{2s+2k-2}-5^{2s-1}-5^{2k-1}+1}{24}\Big)q^n\\
		=\sum_{i=1}^{\infty}\tau_{2\alpha+3,i}\,q^{i-1}\Big(\frac{E_5}{E_1}\Big)^{6i-2}
	\end{align*}
	which is the $\alpha+1$ case of (\ref{gf31}).\\This completes the proof of Theorem \ref{3.1} by induction.
\end{proof}
	%%%%%%%%%%%%%%%%%%%%%%%%%%%%%%%%%%%%%%%%%%%%%%%%%%%%%%%%%%%%%%%%%%%%%%%%%%%%%%%%%%

	%%%%%%%%%%%%%%%%%%%%%%%%%%%%%%%%%%%%%%%%%%%%%%%%%%%%%%%%%%%%%%%%%%%%%%%%%%%%%%%
	\section{Proof of Theorems \ref{T1} - \ref{t3} }
	We are now ready to prove Theorems 1.2 to 1.4.
	For any integer $n$, let $v(n)$ be the largest power of $5$ that divides $n$. Note that $v(0)=\infty$.
	
	\begin{lemma}[{\cite[Eq.~(6.5.1)]{power}}]\label{4.1}
		For any $i\geq 1$ and $j\geq 1$,
		\begin{equation}\label{2.10}
			v(m_{i,j}) \geq \left[\frac{5j-i-1}{2}\right].
		\end{equation}
	\end{lemma}

	%%%%%%%%%%%%%%%%%%%%%%%%%%%%%%%%%%%%%%%%%%%%%%%%%%%%%%%%%%%%%%%%%%%%%%%%%
	\begin{lemma}\label{4.2}
		For any $k\geq 1$, $0\leq\alpha\leq{s-k-1}$ and $j\geq 1$,
		\begin{equation}\label{yvo}
			v\big(y^{(2k-1)}_{2\alpha+1,j}\big) \geq 2k+\alpha-1+\left[\frac{5j-5}{2}\right]
		\end{equation}
		and
		\begin{equation}\label{yve}
			v\big(y^{(2k-1)}_{2\alpha+2,j}\big) \geq 2k+\alpha-1+\delta_{1,j}+\left[\frac{5j-6}{2}\right],
		\end{equation}
		where
		\begin{align*}
			\delta_{1,j}=\begin{cases}
				1 \qquad\text{if $j=1$},\\
				0 \qquad \text{if $j\neq1$}.
			\end{cases}
		\end{align*}
	\end{lemma}
	
	\begin{proof}
		It follows from \cite[Eq.~(4.3)]{Hunt} that
		\[
		v\big(y^{(2k-1)}_{1,j}\big)=v(x_{2k-1,j}) \geq 2k-1+\left[\frac{5j-5}{2}\right].
		\]
		This implies that \eqref{yvo} holds for $\alpha=0$.
		
		Now suppose \eqref{yvo} is true for some $\alpha\geq0$. Then  by (\ref{Ydef}), we obtain that
		\begin{align*}
			v\big(y^{(2k-1)}_{2\alpha+2,j}\big)
			&\geq \min_{i\geq 1}\left\{v\big(y^{(2k-1)}_{2\alpha+1,i}\big)+v(a_{i,j})\right\}\\
			&\geq \min_{i\geq 1}\left\{2k+\alpha-1+\left[\frac{5i-5}{2}\right]+\left[\frac{5j-i-5}{2}\right]+\delta_{1,i}\,\delta_{1,j}\right\}\\
			&\geq 2k+\alpha-1+\delta_{1,j}+\left[\frac{5j-6}{2}\right],
		\end{align*}
		which is \eqref{yve}.
		
		Now assume that \eqref{yve} holds for some $\alpha\geq 0$. Then,
		\begin{align*}
			v\big(y^{(2k-1)}_{2\alpha+3,j}\big)
			&\geq \min_{i\geq 1}\left\{v\big(y^{(2k-1)}_{2\alpha+2,i}\big)+v(b_{i,j})\right\}\\
			&\geq \min_{i\geq 1}\left\{2k+\alpha-1+\delta_{1,i}+\left[\frac{5i-6}{2}\right]+\left[\frac{5j-i-1}{2}\right]\right\}\\
			&\geq 2k+\alpha+\left[\frac{5j-5}{2}\right].
		\end{align*}
		This implies \eqref{yvo} holds for $\alpha+1$. We therefore complete the proof by induction.
	\end{proof}
	%%%%%%%%%%%%%%%%%%%%%%%%%%%%%%%%%%%%%%%%%%%%%%%%%%%%%%%%%%%%%%%%%%%%%%%%%%%%%%%%
\begin{comment}
		\begin{proof}
		It follows from \cite[Eq.~(4.3)]{Hunt} that
		\[
		v\big(y^{(2k-1)}_{1,j}\big)=v(x_{2k-1,j}) \geq 2k-1+\left[\frac{5j-5}{2}\right].
		\]
		This implies that \eqref{yvo} holds for $\alpha=0$.
		
		Now we assume that \eqref{yvo} is true for some $\alpha>0$. With the help of (\ref{Ydef}), we obtain that
		\begin{align*}
			v\big(y^{(2k-1)}_{2\alpha+2,j}\big)
			&\geq \min_{i\geq 1}\left\{v\big(y^{(2k-1)}_{2\alpha+1,i}\big)+v(a_{i,j})\right\}\\
			&\geq \min_{i\geq 1}\left\{k+\alpha+\left[\frac{5i-5}{2}\right]+\left[\frac{5j-i-2}{2}\right]\right\}\\
			&\geq k+\alpha+1+\left[\frac{5j-5}{2}\right],
		\end{align*}
		which is \eqref{2.12}.
		
		Now suppose \eqref{2.12} holds for some $\alpha\geq 0$. Then,
		\begin{align*}
			v\big(y^{(2k-1)}_{2\alpha+3,j}\big)
			&\geq \min_{i\geq 1}\left\{v\big(y^{(2k-1)}_{2\alpha+2,i}\big)+v(a_{i,j})\right\}\\
			&\geq \min_{i\geq 1}\left\{k+\alpha+1+\left[\frac{5i-5}{2}\right]+\left[\frac{5j-i-2}{2}\right]\right\}\\
			&\geq k+\alpha+1+\left[\frac{5j-5}{2}\right].
		\end{align*}
		This implies \eqref{2.11} holds for $\alpha+1$. We therefore complete the proof by induction.
	\end{proof}
\end{comment}
	%%%%%%%%%%%%%%%%%%%%%%%%%%%%%%%%%%%%%%%%%%%%%%%%%%%%%
	
	\begin{lemma}\label{4.3}
		For any $s>k\geq1$, $1\leq\alpha\leq 2k$ and $j\geq1$,
		\begin{equation}\label{zv}
			v\big(z_{\alpha,j}\big) \geq k+s-1+\left[\frac{\alpha-1}{2}\right]+\left[\frac{5j-5}{2}\right].
		\end{equation}
	\end{lemma}

	%%%%%%%%%%%%%%%%%%%%%%%%%%%%%%%%%%%%%%%%%%%%%%%%%%%%%%
\begin{proof}
	From Lemma \ref{4.2} we have
	\[
	v\big(y^{(2k-1)}_{2s-2k,i}\big) \geq k+s-2+\delta_{1,i}+\left[\frac{5i-6}{2}\right].
	\]
	Since $Z_1=Y^{(2k-1)}_{2s-2k}B$, we obtain that
	\begin{align*}
		v\big(z_{1,j}\big)
		&\geq \min_{i\geq 1}\left\{v\big(y^{(2k-1)}_{2s-2k,i}\big)+v(b_{i,j})\right\}\\
		&\geq \min_{i\geq 1}\left\{k+s-2+\delta_{1,i}+\left[\frac{5i-6}{2}\right]+\left[\frac{5j-i-1}{2}\right]\right\}\\
		&\geq k+s-1+\left[\frac{5j-5}{2}\right].
	\end{align*}
	This implies that \eqref{zv} holds for $\alpha=1$.
	
	Now suppose \eqref{zv} is true for some odd $\alpha$ with $1\leq\alpha<2k$. Then by (\ref{zdef}), we obtain that
	\begin{align*}
		v\big(z_{\alpha+1,j}\big)
		&\geq \min_{i\geq 1}\left\{v\big(z_{\alpha,i}\big)+v(c_{i,j})\right\}\\
		&\geq \min_{i\geq 1}\left\{k+s-1+\left[\frac{\alpha-1}{2}\right]+\left[\frac{5i-5}{2}\right]+\left[\frac{5j-i-4}{2}\right]\right\}\\
		&\geq k+s-1+\left[\frac{\alpha}{2}\right]+\left[\frac{5j-5}{2}\right],
	\end{align*}
	which is \eqref{zv} for $\alpha+1$.
	
	Now assume that \eqref{zv} holds for some even $\alpha$ with $2\leq\alpha<2k$. Then,
	\begin{align*}
		v\big(z_{\alpha+1,j}\big)
		&\geq \min_{i\geq 1}\left\{v\big(z_{\alpha,i}\big)+v(b_{i,j})\right\}\\
		&\geq \min_{i\geq 1}\left\{k+s-1+\left[\frac{\alpha-1}{2}\right]+\left[\frac{5i-5}{2}\right]+\left[\frac{5j-i-1}{2}\right]\right\}\\
		&\geq k+s-1+\left[\frac{\alpha}{2}\right]+\left[\frac{5j-5}{2}\right].
	\end{align*}
	This implies \eqref{zv} holds for $\alpha+1$. We therefore complete the proof by induction.
\end{proof}
	\begin{lemma}\label{4.4}
			For any $s>k\geq 1$, $\alpha\geq 1$ and $j\geq 1$,
		\begin{equation}\label{2.12}
		v\big(\tau_{\alpha,j}\big) \geq 2k+s-2+\left[\frac{\alpha-1}{2}\right]+\left[\frac{5j-5}{2}\right].
		\end{equation}
	\end{lemma}
		\begin{proof}
			From Lemma \ref{4.3} we have
			\[
			v\big(\tau_{1,j}\big)=v\big(z_{2k-1,j}\big) \geq 2k+s-2+\left[\frac{5j-5}{2}\right].
			\]
			This implies that \eqref{2.12} holds for $\alpha=1$.
			
			Now we assume that \eqref{2.12} is true for some odd $\alpha\geq1$. Then by (\ref{tdef}), we obtain that
			\begin{align*}
				v\big(\tau_{\alpha+1,j}\big)
				&\geq \min_{i\geq 1}\left\{v\big(\tau_{\alpha,i}\big)+v(c_{i,j})\right\}\\
				&\geq \min_{i\geq 1}\left\{2k+s-2+\left[\frac{\alpha-1}{2}\right]+\left[\frac{5i-5}{2}\right]+\left[\frac{5j-i-4}{2}\right]\right\}\\
				&\geq 2k+s-2+\left[\frac{\alpha}{2}\right]+\left[\frac{5j-5}{2}\right],
			\end{align*}
			which is \eqref{2.12} for $\alpha+1$.
			
			Now suppose \eqref{2.12} holds for some even $\alpha\geq2$. Then,
			\begin{align*}
				v\big(\tau_{\alpha+1,j}\big)
				&\geq \min_{i\geq 1}\left\{v\big(\tau_{\alpha,i}\big)+v(d_{i,j})\right\}\\
				&\geq \min_{i\geq 1}\left\{2k+s-2+\left[\frac{\alpha-1}{2}\right]+\left[\frac{5i-5}{2}\right]+\left[\frac{5j-i-2}{2}\right]\right\}\\
				&\geq 2k+s-2+\left[\frac{\alpha}{2}\right]+\left[\frac{5j-5}{2}\right].
			\end{align*}
			This implies \eqref{2.12} holds for $\alpha+1$. We therefore complete the proof by induction.
		\end{proof}
	\begin{proof}[\textbf{Proof of Theorem \ref{T1}}]
		It follows from (\ref{yvo}) and (\ref{yve}) that
	\begin{equation}\label{4.6}
		v\big(y^{(2k-1)}_{\beta,j}\big)\ \ge\
		\begin{cases}
			2k+\alpha-1, & j=1,\\[2pt]
			2k+\alpha+1, & j\ge 2,
		\end{cases}
		\qquad \beta\in\{2\alpha+1,\,2\alpha+2\}.
	\end{equation}
		The congruence (\ref{t11}) follows immediately from above inequality and (\ref{gf1}).
		
		Moreover, we deduce that, modulo $5^{2k+\alpha}$,
		\begin{align}
				\sum_{n=0}^{\infty}RD_{(5^{2k-1},5^{2s-1})}\big(5^{2k+2\alpha-1}n+\frac{4\cdot5^{2k+2\alpha}-5^{2k-1}+1}{24}\big)q^n 
			\equiv\frac{{E_{5^{2(s-k)-2\alpha}}}}{{E_{5^{2s-2\alpha-1}}}}\cdot y_{2\alpha+1,1}^{2k-1}\cdot E_5^4 .\label{extra11}
		\end{align}
		Since no powers of $q$ other than 5 are present,	the congruences \eqref{t12} follows.
		This completes the proof of Theorem \ref{T1}.
	\end{proof}
	\begin{proof}[\textbf{Proof of Theorem \ref{t2}}]
		It follows from Lemma \ref{4.3} that
		\begin{align}
			\nu\big(z_{\alpha,j}\big) &\geq \label{vrqz1}
			\begin{cases}
				k+s-1+\left[\frac{\alpha-1}{2}\right] & \text{if } j=1,\\
				k+s+1+ \left[\frac{\alpha-1}{2}\right]& \text{if } j\geq 2.
			\end{cases} 
		\end{align}
		The congruence \ref{t21} follows immediately from above inequality and (\ref{gf21}).
		
		Moreover, we deduce that, modulo $5^{k+s+\alpha}$,
		\begin{align}
		\sum_{n=0}^{\infty}RD_{(5^{2k-1},5^{2s-1})}\Big(5^{2s+2\alpha-1}n+\frac{21\cdot5^{2s+2\alpha-1}-5^{2s-1}-5^{2k-1}+1}{24}\Big)q^n\nonumber\\
		\equiv\frac{1}{E_{5^{2k-2\alpha-1}}}z_{2\alpha+1,1}{E_5^{4}}{E_1}.\label{extra1}
		\end{align}
		Thanks to (\ref{r(q)})
		\begin{equation}\label{2.24}
			E_1=E_{25}\big(R(q^{5})^{-1}-q-q^{2}R(q^{5})\big),
		\end{equation}
		where
		\[
		R(q)=\frac{(q;q^{5})_{\infty}(q^{4};q^{5})_{\infty}}{(q^{2};q^{5})_{\infty}(q^{3};q^{5})_{\infty}}.
		\]
		The congruences (\ref{t23}) and (\ref{t24}) follow by combining (\ref{extra1}) and \eqref{2.24}.\\
		In a similar manner congruence (\ref{t22}) follows from (\ref{gf22}) and (\ref{vrqz1}).
		This completes the proof of Theorem \ref{t2}.
	\end{proof}
		\begin{proof}[\textbf{Proof of Theorem \ref{t3}}]
		It follows from Lemma \ref{4.4} that
		\begin{align}
			v\big(\tau_{\alpha,j}\big) &\geq \label{vrqt1}
			\begin{cases}
				2k+s-2+\left[\frac{\alpha-1}{2}\right] & \text{if } j=1,\\[2pt]
				2k+s+\left[\frac{\alpha-1}{2}\right] & \text{if } j\geq 2.
			\end{cases}
		\end{align}
		The congruence (\ref{t31}) follows immediately from above inequality and (\ref{gf31}).
		
		Moreover, we deduce that, modulo $5^{2k+s+\alpha-1}$,
		\begin{align}
			\sum_{n=0}^{\infty}RD_{(5^{2k-1},5^{2s-1})}\Big(5^{2s+2k+2\alpha-3}n+\frac{16\cdot5^{2s+2k+2\alpha-3}+5^{2s+2k-2}-5^{2s-1}-5^{2k-1}+1}{24}\Big)q^n\nonumber\\
			\equiv\tau_{2\alpha+1,1}\,{E_5^{3}}\,{E_1}.\label{extra31}
		\end{align}
		Thanks to (\ref{r(q)})
		\begin{equation}\label{2.24b}
			E_1=E_{25}\big(R(q^{5})^{-1}-q-q^{2}R(q^{5})\big),
		\end{equation}
		where
		\[
		R(q)=\frac{(q;q^{5})_{\infty}(q^{4};q^{5})_{\infty}}{(q^{2};q^{5})_{\infty}(q^{3};q^{5})_{\infty}}.
		\]
		The congruences (\ref{t33}) follow by combining (\ref{extra31}) and \eqref{2.24b}.\\
		In a similar manner congruence (\ref{t32}) follows from (\ref{gf32}) and (\ref{vrqt1}). Moreover, we deduce that, modulo $5^{2k+s+\alpha-1}$,
		\begin{align}
			\sum_{n=0}^{\infty}RD_{(5^{2k-1},5^{2s-1})}\Big(5^{2s+2k+2\alpha-2}n+\frac{8\cdot5^{2s+2k+2\alpha-2}+5^{2s+2k-2}-5^{2s-1}-5^{2k-1}+1}{24}\Big)q^n\nonumber\\
			\equiv\tau_{2\alpha+2,1}\,{E_5}\,{E_1^{3}}.\label{extra32}
		\end{align}
			With the help of \eqref{r(q)}, we obtain that
		\begin{equation}\label{f1^3}
			E_1^{3}=E_{25}^{3}\big(R(q^{5})^{-3}-3qR(q^{5})^{-2}+5q^{3}-3q^{5}R(q^{5})^{2}-q^{6}R(q^{5})^{3}\big).
		\end{equation}
		the congruence (\ref{t34}) follow by combining (\ref{extra32}) and \eqref{f1^3}.
		This completes the proof of Theorem \ref{t3}.
	\end{proof}

	\noindent\textsuperscript{1,2}Ramanujan School of Mathematical Sciences,\\Department of Mathematics,\\ Pondicherry University,\\ Puducherry - 605014, India. \\
	
	\noindent Email: \texttt{anjelinvallialil@pondiuni.ac.in} \\
	Email: \texttt{dr.fathima.sn@pondiuni.ac.in}  (\Letter)
\end{document}